\documentclass[a4paper,11pt]{article}

\usepackage{amssymb,latexsym,amsmath,amsthm,amsfonts,graphics}

\def\B{\mathcal{B}}
\def\R{\mathbb{R}}
\def\N{\mathbb{N}}
\def\Z{\mathbb{Z}}
\def\E{\mathbb{E}}
\def\P{\mathbb{P}}

\def\S{\mathbb{S}}
\def\C{\mathbb{C}}

\def\HH{\mathcal{H}}

\def\XX{\tilde{X}}

\def\a{\vec{a}}
\def\b{\vec{b}}
\def\c{\vec{c}}
\def\e{\vec{e}}

\def\i{\mathbf{i}}

\def\prf{\vskip -0.1cm\noindent{\bf Proof\quad }}
\def\prfend{\hfill{$\Box$}\vskip 0.2cm}
\newtheorem{lem}{Lemma}

\newtheorem{thm}{Theorem}

\newtheorem{exa}[lem]{Example}

\date{}
\title{Width deviation of convex polygons\\
\author{Shigeki Akiyama\thanks{Institute of Mathematics, University of Tsukuba (akiyama@math.tsukuba.ac.jp)}\quad and Teturo Kamae\thanks{Advanced Mathematical Institute, Osaka City University (kamae@apost.plala.or.jp)}}}

\begin{document}
\maketitle
\begin{abstract}
We consider the 
width $X_T(\omega)$ of a convex $n$-gon $T$ in the plane along the random direction $\omega\in\mathbb{R}/2\pi \mathbb{Z}$
and study its deviation rate:
$$
\delta(X_T)=\frac{\sqrt{\E(X^2_T)-\E(X_T)^2}}{\E(X_T)}.
$$
We prove that the maximum is attained if and only if 
$T$ degenerates to a $2$-gon.
Let $n\geq 2$ be an integer which is not a power of $2$. 
We show that $$
\sqrt{\frac{\pi}{4n\tan(\frac{\pi}{2n})}
+\frac{\pi^2}{8n^2\sin^2(\frac{\pi}{2n})}-1}
$$
is the minimum of $\delta(X_T)$ among all $n$-gons and determine completely the shapes of $T$'s which attain this minimum.
They are characterized as polygonal approximations 
of equi-Reuleaux bodies, found and studied by K.~Reinhardt \cite{RH}.
In particular, if $n$ is odd, then the
regular $n$-gon is one of the minimum shapes. 
When $n$ is even, we see that regular $n$-gon is far from optimal.
We also observe an unexpected property of the deviation rate on the truncation of the regular triangle.
\end{abstract}

\section{Introduction}
The width of an image along various directions is basic information in the 
Image Processing Technique. We are interested in the deviation of the 
widths of compact convex sets in the plane along a random direction. 
Although it is an important and useful quantity both practically and theoretically, 
as long as we know, there is no serious study on it. 

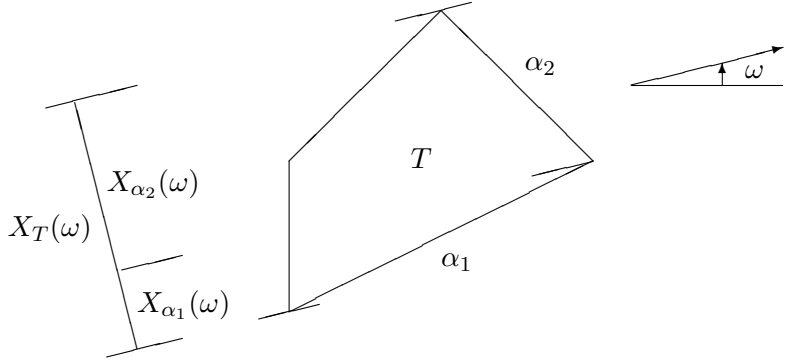
\begin{figure}[h]
\setlength{\unitlength}{1mm}
\begin{picture}(80,50)(-67,-45)
\put(20,-20){\line(-1,1){20}}
\put(0,0){\line(-1,-1){20}}
\put(-20,-20){\line(0,-1){20}}
\put(-20,-40){\line(2,1){40}}
\put(4,1){\line(-4,-1){10}}
\put(-52,-13){\line(4,1){12}}

\put(20,-20){\line(-4,-1){8}}
\put(-42,-34.5){\line(4,1){8}}
\put(-40,-40){$X_{\alpha_1}(\omega)$}
\put(-44,-24){$X_{\alpha_2}(\omega)$}
\put(0,-34){$\alpha_1$}
\put(11,-8){$\alpha_2$}

\put(-16,-39){\line(-4,-1){8}}
\put(-44,-46){\line(4,1){10}}
\put(-40,-45){\line(-1,4){8.2}}
\put(-4,-21){$T$}
\put(-57,-30){$X_T(\omega)$}
\put(40,-9){$\omega$}
\put(37,-10){\vector(0,1){3}}
\put(25,-10){\vector(4,1){20}}
\put(25,-10){\line(1,0){20}}
\end{picture}
\caption{$X_T(\omega)$}
\end{figure}

Let $\Omega=\R/2\pi\Z=(-\pi,\pi]$ and $\P$ be the normalized 
Lebesgue measure on $\Omega$. 
Consider the probability space $(\Omega,\P)$; we write $\Omega$ for short.
Let a compact convex set $T$ in the complex plane $\C$ be given. 
We identify $T$ with its boundary. 
We discuss the length of the orthogonally projected shadow of $T$ 
by the light from a random direction $\omega$, say $X_T(\omega)$. 
It is also interpreted as the width of $T$ along the orthogonal direction of $\omega$, 
that is, $\omega\pm\frac{\pi}{2}$. 

We are interested in the uniformity (or its contrary) of the deviation of 
$X_T(\omega)$ with respect to $\omega\in\Omega$. 
For this purpose, we consider the {\it deviation rate} of  the random variable $X_T$, 
say $\delta(X_T)$, defined by
\begin{equation}
\delta(X_T)=\frac{\sqrt{\E(X^2_T)-\E(X_T)^2}}{\E(X_T)},
\end{equation}
where $\E(~)$ is the expectation of the random variables. 
Clearly, $\delta(X_T)$ is invariant among the similarity images of $T$. 
It is also clear that $\delta(X_T)=0$ if and only if $T$ is a closed curve of constant 
width (c.f. \cite{MMO}). Hence, $\delta(X_T)$ measures the degree how far $T$ is from convex bodies of constant width. 
For the ellipse 
$T=\{(x,y); ~(x^2/a^2)+(y^2/b^2)=1\}$ with the perimeter $L$, it is not difficult 
to see that 
$$\delta(X_T)=\frac{\sqrt{2\pi^2(a^2+b^2)-L^2}}{L}.$$

In this paper, we study the deviation rate of the 
convex $n$-gons $T$ with $n\ge 2$. We decompose $X_T$ as the sum 
of random variables coming from its edges, say 
$X_T=\sum_{j=1}^nX_{\alpha_j}$. 
Throughout this paper, the branch of the argument of a complex number is chosen to be $(-\pi,\pi]$.
For a complex number $\alpha\in\C$ with $\arg(\alpha)=\theta$, define a random 
variable $X_\alpha$ on $\Omega$ by
$$ 
X_\alpha(\omega)=|\alpha|\sin(\theta-\omega)_+,
$$
where $x_+=\max\{x, 0\}$. In other words,
$X_\alpha(\omega)$ is the length of the orthogonally 
projected shadow of the right side of the vector $\vec{O\alpha}$ 
by the light from the $\omega$ direction (the left side of the vector is permeable and 
makes no shadow). 

Let $\alpha_1,\alpha_2,\cdots,\alpha_n~(n\ge 2)$ be a sequence of distinct complex 
numbers. They are arranged in the counter clockwise order forming a convex $n$-gon 
if and only if 
\begin{multline*}
(*)~\arg(\alpha_2-\alpha_1)\le\arg(\alpha_3-\alpha_2)\le\cdots\\
\le\arg(\alpha_{n+1}-\alpha_n)\le\arg(\alpha_2-\alpha_1)+2\pi~~({\rm mod}~2\pi),
\end{multline*}
where \underline{we always consider the suffix $j$ related to the $n$-gon $T$ 
in modulo $n$}, so that $\alpha_{n+1}=\alpha_1,~\alpha_0=\alpha_n$, etc.
In this case, the convex $n$-gon with vertices $\alpha_1,\cdots,\alpha_n$ 
is denoted by $T=T(\alpha_1,\cdots,\alpha_n)$. 

The above $T$ is non-degenerate (i.e. all 
the vertices are extremal points) if and only if ``$<$'' holds everywhere 
in the above inequalities. If this is not the case, then we identify 
$T=T(\alpha_1,\cdots,\alpha_n)$ with $T(\alpha'_1,\cdots,\alpha'_m)$, 
where $\{\alpha'_1,\cdots,\alpha'_m\}$
are the set of extremal points among $\{\alpha_1,\cdots,\alpha_n\}$ arranged in the 
counter clockwise order.

Define a random variable $X_T$ by
$$
X_T(\omega)=\sum_{j=1}^nX_{\alpha_{j+1}-\alpha_j}(\omega),
$$
which agrees with the former explanation as to the length of the shadow of $T$ (see Figure 1). 
  
The set of convex $n$-gons $T=T(\alpha_1,\cdots,\alpha_n)$ can be identified 
with the sequence of $n$ distinct complex numbers $(\alpha_1,\cdots,\alpha_n)$ 
satisfying $(*)$. We denote the space consisting of these 
$(\alpha_1,\cdots,\alpha_n)$ by $\Theta_n$. We consider the 
usual structures coming from $\C^n$ on $\Theta_n$. 
Then, $\cup_{k=2}^{n-1}\Theta_k$ is considered as the boundary of $\Theta_n$ 
in the above sense. 

Denote the factor space of $\Theta_n$ divided by the similarity equivalence 
by $\tilde{\Theta}_n$. Then it is a compact space. 
Most of our notions like the regular $n$-gon are notions in $\tilde{\Theta}_n$ 
rather than $\Theta_n$. The deviation rate $\delta(X_T)$ 
can be considered as a continuous  and piecewise smooth functional on $\tilde{\Theta}_n$. 
Hence, it has the minimum and the maximum in $\Theta_n$. 

It turned out that
the solution of the minimization problem of $\delta(X_T)$ is related to
the convex bodies of constant width.
A {\it Reuleaux body}\footnote{In literature, it is often referred to Reuleaux "polygon", 
though its boundary is not linear. 
In this paper, we use the word "body" to distinguish it from a genuine polygon.} is 
a convex body of constant width whose boundary consists of a finite number of circular arcs with the center in it and the radius equal to the width. 

The reader can consult \cite{MMO, Mo}, but for 
the self-containedness, we review the Reuleaux body here briefly. 

Let $D$ be a Reuleaux body with width $r$. Then, $C=\partial D$ consists of 
a finite number of circular arcs, say $C_1,\cdots,C_p$, of radius $r$. 
The endpoints of the circular arcs are called {\it essential} if the 2 neighboring circular 
arcs have different centers, and hence cannot be joined into a single circular arc. 

We assume that the circular arcs $C_1,\cdots,C_p$ are arranged in the 
counter clockwise order and all the endpoints are essential. 
Let $A_1,A_2$ be the endpoints of $C_1$; $A_2,A_3$ be the endpoints of $C_2$, 
$\cdots$; $A_p,A_{p+1}$ be the endpoints of $C_p$ (the suffices of $C$ or $A$ are  
considered modulo $p$ so that $A_{p+1}=A_1$). 
Note that a center of any one of
the circular arcs is an essential endpoint of some other circular arcs on $C$, 
since otherwise, we have 2 points in $C$ with a distance larger than $r$. 
Conversely 
any endpoint is a center of some circular arc since 
there must be the same number of circular arcs and centers. 

Hence, the center of $C_1$ is one of $A_3,\cdots,A_p$. Let it be $A_k$. 
Since $D$ is strictly convex, the tangent vector of $C$ at $c\in C$ 
to the counter clockwise direction rotates counter clockwise as $c$ moves. 
The center of the circular arc containing $c\in C$ is the other intersection of the 
normal line to the tangent vector at $c$ with $C$. This intersection moves at the endpoints of the circular arcs to the counter clockwise direction since the normal lines from the same point rotate to the counter clockwise direction. Hence, the center 
of this circular arc is the next end point to that before. 
This means that the center of $C_1$ is $A_k$, the center of $C_2$ is $A_{k+1}$, 
$\cdots$. 

On the other hand, the center of the circular arc $C_k$ must be $A_2$ since 
otherwise, either it is $A_1$ or there is $l\ne1,2$ such that 
$\overline{A_kA_l}=r$. The latter is impossible since if so, by the convexity and the assumption on $C$, the circular arc $C_1$ can be extended to $A_l$ contradicting 
the assumption that both $A_1$ and $A_2$ are essential end points. The former is impossible since if so, then we have a contradiction that $\overline{A_2A_{k+1}}>r$. 
Thus, we have 
$$
A_2=\mbox{the center of }C_k=A_{2k-1},
$$
and hence, $2\equiv 2k-1~(\mbox{mod}~p)$. This implies that $p$ is odd. 

Consider the diagonals connecting 2 points in $C$, one of the center and 
an interior point of the circular arc centered by it. Any of 2 diagonals with different 
centers always intersect just at one point. 
Let $\Lambda_j$ be the set of angles from $A_j$ to the points in the circular arc 
centered by $A_j$. We always consider the angles in modulo $2\pi$. 
Then, by the above argument, $\Lambda_j$ and $\Lambda_l$ with $j\ne l$ 
are essentially disjoint. The same thing holds for $\Lambda_j$ and $\Lambda_l+\pi$. 
Let $\Lambda=\cup_{j=1}^p\Lambda_j$. Then, the union is essentially disjoint. 
Moreover, $\Lambda$ and $\Lambda+\pi$ is also essentially disjoint. 
It also holds that $\Lambda\cup(\Lambda+\pi)=(-\pi,\pi]$. 
This is because the Lebesgue measure of $\Lambda$ is $\pi$, since 
the integration of the exterior angle of $C$ is $2\pi$, which counts the 
angles in $\Lambda_j~(j=1,\cdots,p)$ twice, once at $C_j$, once at its center, 
where the exterior angle jumps just the amount of angles in $\Lambda_j$. 

Thus, it holds that the circular arcs $C_1,\cdots,C_p$ can be embedded by parallel 
translations by 
$\vec{v}_1,\cdots,\vec{v}_p\in\C$ into a circle $\S$ of radius $r$ so that 
$$
\tilde{\Lambda}:=\cup_{j=1}^p(C_j+\vec{v}_j)\subset\S
$$
satisfies that 
\begin{align}
&\#((C_j+\vec{v}_j)\cap(C_k+\vec{v}_k))<\infty~\mbox{for any}~j\ne k\nonumber,\\
&\#(\tilde{\Lambda}\cap R_\pi\tilde{\Lambda})<\infty
~\mbox{and}~\tilde{\Lambda}\cup R_\pi\tilde{\Lambda}=\S,
\end{align}
where $R_\pi$ is the rotation on $\S$ by angle $\pi$, and ``$\#$'' implies the number 
of elements in a set. 

We write a {\it Reuleaux $p$-body} to indicate the number $p$ of circular 
arcs with different centers. In particular, a Reuleaux $p$-body is called 
{\it regular} if all the circular arcs have the same length. A regular Reuleaux 3-body 
is well known as ``Reuleaux triangle''. 


A {\it Reinhardt $n$-gon} (\cite{HM13}, \cite{RH}) is an equilateral convex 
$n$-gon that can be inscribed in a Reuleaux body, containing all the essential 
endpoints of the circular arcs.  
It is referred to as a {\it Reinhalt polygon} 
if $n$ is not necessarily specified. 
If the lengths of the circular arcs of a Reuleaux $p$-body 
have ratio $n_1:n_2:\dots:n_p$ with integer $n_i$'s, 
we can divide the circular arcs into $n_i\ (i=1,\dots,p)$ parts of equal length
and take the convex-hull of all the division points to get a Reinhardt $n$-gon, 
where $n=n_1+\cdots+n_p$, which is specially called a Reinhardt $(n_1,\cdots,n_p)$-gon. 

\begin{figure}[h]
\begin{center}
\includegraphics[13cm,4cm]{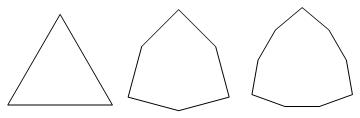}
\end{center}
\vspace{-2em}
\caption{Reinhardt $i^3$-gon for $i=1,2,3$\label{Reu}}
\end{figure}

In particular, let $p\ge 3$ be an odd integer and take a regular Reuleaux 
$p$-body $D$.
We divide all the circular arcs of $D$ into $q$ parts of equal length.
Then, the convex hull of all the division points forms a Reinhardt $pq$-gon, 
which is specially a Reinhardt $(\underbrace{q,\cdots,q}_{\mbox{$p$ times}})$-gon, 
see Figure \ref{Reu}. Denote $q^p=(\underbrace{q,\cdots,q}_{\mbox{$p$ times}})$ for short. 
Note that Reinhardt $1^p$-gon is the regular 
$p$-gon.

A Reinhardt polygon gives the solution of three optimization problems
on convex $n$-gons when $n$ is not a power of $2$:

\begin{description}
    \item[Problem 1.] Maximize the perimeter for a fixed diameter \cite{RH}, 
    \item[Problem 2.] Maximize the width for a fixed diameter \cite{BF},
    \item[Problem 3.] Maximize the width for a fixed perimeter \cite{AHM}.
\end{description}

Given a cyclic integer vector $(n_1,n_2,\dots,n_p)$, there exists a unique way 
to construct a Reinhardt $(n_1,n_2,\dots,n_p)$-gon when this is possible.  
We reproduce a necessary and sufficient condition of Reinhardt \cite{RH}
that a given integer
cyclic vector $(n_1,n_2,\dots,n_p)$ 
forms a Reinhardt $(n_1,n_2,\dots,n_p)$-gon as a byproduct of our result
(Theorem \ref{Main} and \ref{MainGeom}). In particular,
when $n$ is not a power of $2$, we have
a Reinhardt $(n/p)^p$-gon 
for any odd divisor $p>1$. 
Many interesting properties on Reinhardt polygons 
are discussed in \cite{Mo, Mo2, HM13,HM19}.
When $n$ is odd, the regular $n$-gon is one of the Reinhardt polygons but it is not unique when $n$ is a composite. Reinhardt $n$-gon  
is unique if and only if $n=p$ or $2p$ with a prime $p$. 
A Reinhardt polygon may have no symmetry at all (see Figures in \cite{HM13}). 
\\

In this paper, we prove the following results. \vspace{0.5em}\\
(1) The deviation rate for the compact convex sets $T$ attains the maximum 
if and only if $T$ is the 2-gon (Theorem \ref{Max}). \vspace{0.5em}\\
(2) For the integer $n\ge 2$ which is not a power of $2$, 
the minimum of $\delta(X_T)$ for $T\in\Theta_n$ is
\begin{equation}
\label{nu}
\nu_n:=\sqrt{\frac{\pi}{4n\tan(\frac{\pi}{2n})}
+\frac{\pi^2}{8n^2\sin^2(\frac{\pi}{2n})}-1}.
\end{equation}

Theorem \ref{Main} gives a complete characterization of 
the shapes that attain the minimum value $\nu_n$.
As a consequence, we show that the minimum shape is nothing but 
a Reinhardt $n$-gon (Theorem \ref{MainGeom}).

The minimum values for odd $n\ge 3$ is strictly decreasing in $n$.
\vspace{0.5em}\\ 
(3) For even $m\ge 4$, regular $m$-gon is far from the minimum shape. 
Let $n(<m)$ be the odd number such that either 
$n=\frac{m}{2}$ or $n=\frac{m}{2}+1$. If $n=\frac{m}{2}$, then 
$\delta(X_{T_m})=\delta(X_{T_n})$ holds, 
and if $n=\frac{m}{2}+1$, then $\delta(X_{T_m})>\delta(X_{T_n})$ holds, where 
$T_m,~T_n$ are the regular $m$-gon and $n$-gon, respectively 
(Theorem \ref{Reg}).\\

\begin{figure}[h]
\begin{center}
\includegraphics[9cm,3.5cm]{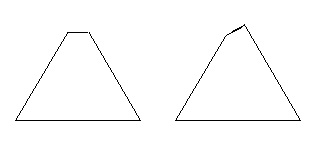}
\end{center}
\vspace{-2em}
\caption{parallel truncation (left) and non-parallel truncation (right)}
\end{figure}

(4) Let $T$ be the regular triangle. We call a {\it truncation} of it a quadrangle obtained 
 by cutting off a vertex by a line near it. It is called a {\it parallel} truncation if 
the line is parallel to the opposite side, otherwise a  {\it non-parallel} truncation 
(Figure 3). We prove that a parallel truncation of the regular triangle increases the deviation rate, while a non-parallel truncation decreases it (Theorem \ref{Trun}).

\section{Maximum of the deviation rate}

\begin{thm}
\label{Fund}
For any $\alpha,\beta\in\C$ with $\eta=\arg(\beta/\alpha)\in(-\pi,\pi]$, 
the following statements hold:\\
$(i)$ $\E(X_\alpha)=|\alpha|/\pi$,\\
$(ii)$ $\E(X_\alpha X_\beta)=\frac{1}{4\pi}|\alpha||\beta|V(\eta)$, 
where $V$ is a periodic function of period $2\pi$ such that 
$$
V(x)=(\pi-|x|)\cos x+\sin|x|~~(-\pi<x\le\pi).
$$
\end{thm}

\begin{figure}[h]
\begin{center}
\includegraphics[6cm,4.5cm]{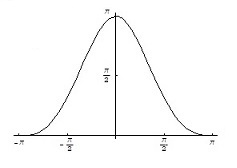}
\end{center}
\vspace{-2em}
\caption{Function $V$}
\end{figure}

\prf
\noindent $(i)$ Let $\theta=\arg(\alpha)$. Then, we have  
$$
\E(X_\alpha)=\frac{1}{2\pi}\int_0^{2\pi}|\alpha|\sin(\theta-\xi)_+d\xi
=\frac{|\alpha|}{2\pi}\int_0^\pi\sin\xi d\xi=\frac{|\alpha|}{\pi}.
$$

\noindent $(ii)$ We have
\begin{align*}
\E(X_\alpha X_\beta)=\frac{1}{2\pi}\int_0^{2\pi}|\alpha||\beta|
\sin(\theta-\xi)_+\sin(\theta+\eta-\xi)_+d\xi
\end{align*}
\begin{align*}
&=\frac{1}{2\pi}|\alpha||\beta|\int_{[0,\pi]\cap[-\eta,\pi-\eta]}
\sin\xi\sin(\xi+\eta)d\xi\\
&=\frac{1}{4\pi}|\alpha||\beta|\int_{[0,\pi]\cap[-\eta,\pi-\eta]}
(\cos\eta-\cos(2\xi+\eta))d\xi.
\end{align*}
Therefore, if $\eta<0$, then
\begin{align*}
&\E(X_\alpha X_\beta)=\frac{1}{4\pi}|\alpha||\beta|\int_{-\eta}^\pi
(\cos\eta-\cos(2\xi+\eta))d\xi\\
&=\frac{1}{4\pi}|\alpha||\beta|((\pi+\eta)\cos\eta-\sin\eta)
=\frac{1}{4\pi}|\alpha||\beta|V(\eta).
\end{align*}
If $\eta\ge 0$, then 
\begin{align*}
&\E(X_\alpha X_\beta)=\frac{1}{4\pi}|\alpha||\beta|\int_0^{\pi-\eta}
(\cos\eta-\cos(2\xi+\eta))d\xi\\
&=\frac{1}{4\pi}|\alpha||\beta|((\pi-\eta)\cos\eta+\sin\eta)
=\frac{1}{4\pi}|\alpha||\beta|V(\eta).
\end{align*}
\prfend

\begin{thm} 
\label{Max}
The maximum of $\delta(X_T)$ for the compact convex sets $T$ 
is attained by a $2$-gon. 
\end{thm}
\prf
Since all 2-gons are similar, they have the same $\delta$-value. 
Let $U=T(0,1)$ be a 2-gon. Then by Theorem \ref{Fund}, we have 
\begin{align*}
&\E(X_U)=\frac{1}{\pi}+\frac{1}{\pi}=\frac{2}{\pi}\\
&\E(X_U^2)=\frac{1}{4\pi}(2V(0)+2V(\pi))=\frac12.
\end{align*}
Hence, $\delta(X_U)=\sqrt{(\pi^2/8)-1}$. 
 
Take any convex polygon $T=T(\alpha_1,\cdots,\alpha_n)$. 
We prove that $\delta(X_T)\le\sqrt{(\pi^2/8)-1}$, and that   
the equality holds only when $T$ is 2-gon. Let 
$$\beta_j=\alpha_{j+1}-\alpha_j=r_je^{\i\theta_j}~(j=1,\cdots,n)$$ 
and let $\theta_{jk}\in(-\pi,\pi]$ satisfies that 
$$\theta_{jk}\equiv\theta_j-\theta_k~({\rm mod}~2\pi)~(j,k=1,\cdots,n).$$

It is sufficient to prove that 
$\frac{\E(X_T^2)}{\E(X_T)^2}\le\frac{\pi^2}{8}$. 
Hence, it is sufficient to 
prove that $I:=\frac{\pi^3}{2}\E(X_T)^2-4\pi\E(X_T^2)\ge0$. 
Then by Theorem \ref{Fund}, we have
\begin{align*}
I=\frac{\pi^3}{2}\E(X_T)^2-4\pi\E(X_T^2)=\frac{\pi^3}{2}
\left(\sum_{j=1}^n\E(X_j)\right)^2-4\pi\sum_{j,k=1}^n\E(X_jX_k)
\end{align*}
\begin{align*}
&=\frac{\pi^3}{2}\left(\sum_{j=1}^n\frac{r_j}{\pi}\right)^2
-\sum_{j,k=1}^nr_jr_kV(\theta_j-\theta_k)=\sum_{j,k=1}^nr_jr_k
\left(\frac{\pi}{2}-V(\theta_j-\theta_k)\right)\\
&=\sum_{j,k=1}^nr_jr_k\left(\frac{\pi}{2}-(\pi-|\theta_{jk}|)\cos\theta_{jk}-|\sin\theta_{jk}|\right).~~~~(\#)
\end{align*}
Since $r_jr_k\cos\theta_{jk}$ is the inner product between 
$\vec{O\beta_j}$ and $\vec{O\beta_k}$, we have
$$
\sum_{j,k=1}^nr_jr_k\cos\theta_{jk}
=\left(\sum_{j=1}^k\vec{O\beta_j},~\sum_{j=1}^k\vec{O\beta_j}\right)
=(\vec{0},\vec{0})=0.
$$
This implies that in the above equality $(\#)$, the constant in the coefficient 
of $\cos\theta_{jk}$ can be changed anyway keeping the equality. Hence, we have
$$
I=\sum_{j,k=1}^nr_jr_k\left(\frac{\pi}{2}-\left(\frac{\pi}{2}-|\theta_{jk}|\right)
\cos\theta_{jk}-|\sin\theta_{jk}|\right).
$$
Thus, to prove $I\ge0$, it is sufficient to prove that 
$$
\left(\frac{\pi}{2}-|x|\right)\cos x+|\sin x|\le\frac{\pi}{2}
$$
for any $x\in(-\pi,\pi]$, which can be verified easily. 
Moreover, the equality holds if and only if $x=0$ or $\pi$, which implies 
that $T$ is 2-gon. 

For a general compact convex set $S$ with $C=\partial S$, 
replacing the above $I$ by  
$$
I=\int\int_{C\times C}
\left(\frac{\pi}{2}-\left(\frac{\pi}{2}-\theta(u,v)\right)\cos\theta(u,v)-|\sin\theta(u,v)|\right)dudv,
$$
where $\theta(u,v)$ is the angle between the tangent lines of $C$ at $u$ and $v$, 
we have the same statement, which completes the proof.   
\prfend

\section{Symmetrization and asymmetrization of $n$-gons}

Let 
$$\HH=\{\alpha\in\C;~\Im(\alpha)>0,~\mbox{or,}~\Im(\alpha)=0~\mbox{and}~\Re(\alpha)\ge0\}$$
be the upper half plane endowed with the quotient topology of $\C$ by identifying $z$ and $-z$. For $\alpha\in\C$, we define 
$$
\iota(\alpha)=\left\{\begin{array}{cc}\alpha&(\alpha\in\HH)\\-\alpha&(\alpha\notin\HH)
\end{array}\right..
$$
For a finite set of nonzero complex numbers $S=\{\alpha_1,\cdots,\alpha_n\}$, 
define 
$\iota(S)\subset\HH$ by 
$$
\iota(S)=\mbox{Abbreviation}\{\iota(\alpha_1),\cdots,\iota(\alpha_n)\},
$$
where $\mbox{Abbreviation}\{\alpha_1',\cdots,\alpha_n'\}$ is the set of complex numbers 
obtained by replacing any of $\alpha_j',\alpha_k'$ with $\arg(\alpha_j')=\arg(\alpha_k')$ 
by $\alpha_j'+\alpha_k'$. 

A nonzero element in $\HH$ is sometimes called a {\it pre-edge}. 
A sequence of pre-edges $\beta_1,\cdots,\beta_m\in\HH$, say 
$\B=\B(\beta_1,\cdots,\beta_m)$ is called a {\it pre-edge bundle} 
(of size $m$) if 
$$
0\le\arg(\beta_1)<\arg(\beta_2)<\cdots<\arg(\beta_m)<\pi.
$$ 
Denote by $\Xi_m$ the set of pre-edge bundles of size $m$. 

For a convex $n$-gon $T=T(\alpha_1,\cdots,\alpha_n)$, we define 
its {\it asymmetrization} $\iota(T)$ as the pre-edge bundle  
$\B=\B(\beta_1,\cdots,\beta_m)$ such that
$$
\iota(\{\alpha_2-\alpha_1,~\alpha_3-\alpha_2,~\cdots,~\alpha_{n+1}-\alpha_n\})
=\{\beta_1,\beta_2,\cdots,\beta_m\}.
$$
In this case, $T$ is called a {\it realization} of $\B$. 
Let $T=T(\alpha_1,\cdots,\alpha_n)$ be a convex $n$-gon and 
$\B=\B(\beta_1,\cdots,\beta_m)$ be its asymmetrization. Then,
\begin{align*}
&U=\\
&T\left(\gamma,\gamma+\frac12\beta_1,\cdots,
\gamma+\frac12(\beta_1+\cdots+\beta_m),
\gamma+\frac12(\beta_2+\cdots+\beta_m),\cdots,\gamma+\frac12\beta_m\right)
\end{align*}
is another realization of $\B$ than $T$, where 
$\gamma=-\frac14(\beta_1+\cdots+\beta_m)$. 
We call $U$ the {\it symmetrization} of $T$ (or $\B$). 

In this case, $U$ is {\it symmetric}, that is, $U=T(\gamma_1,\cdots,\gamma_{2k})$ 
with even size $2k$ and 
$$
\gamma_{k+1}=-\gamma_1,~\gamma_{k+2}=-\gamma_2,~\cdots,
~\gamma_{2k}=-\gamma_k
$$
holds. For a symmetric $U=T(\gamma_1,\cdots,\gamma_{2k})$, it holds that  
\begin{multline}
X_U(\omega)=2|\gamma_j|\sin(\arg(\gamma_j)-\omega)\\
\mbox{for any}~\omega~\mbox{with}
~\arg(\gamma_{j-1}-\gamma_j)\le\omega\le\arg(\gamma_j-\gamma_{j+1}).
\end{multline}
This representation of $X_U$ is called the {\it diagonal representation}. 

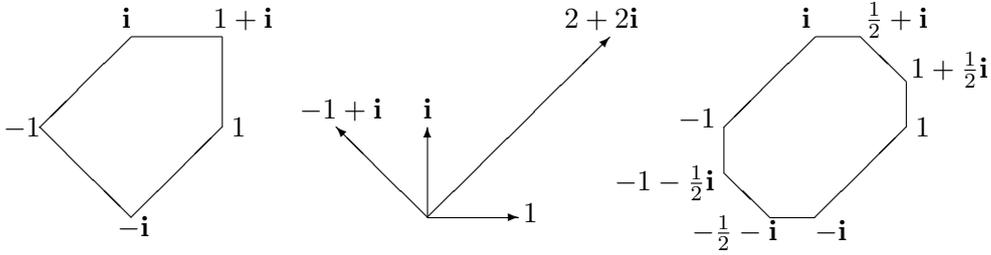
\begin{figure}[h]
\setlength{\unitlength}{0.6mm}
\begin{picture}(200,60)(0,-5)
\put(0,20){\line(1,-1){20}}
\put(20,0){\line(1,1){20}}
\put(40,20){\line(0,1){20}}
\put(40,40){\line(-1,0){20}}
\put(20,40){\line(-1,-1){20}}

\put(-8,18){$-1$}
\put(17,-4){$-\i$}
\put(42,18){$1$}
\put(38,42){$1+\i$}
\put(18,42){$\i$}

\put(85,0){\vector(1,0){20}}
\put(85,0){\vector(1,1){40}}
\put(85,0){\vector(0,1){20}}
\put(85,0){\vector(-1,1){20}}

\put(106,-1){$1$}
\put(115,42){$2+2\i$}
\put(84,22){$\i$}
\put(57,22){$-1+\i$}

\put(160,0){\line(1,0){10}}
\put(170,0){\line(1,1){20}}
\put(190,20){\line(0,1){10}}
\put(190,30){\line(-1,1){10}}
\put(180,40){\line(-1,0){10}}
\put(170,40){\line(-1,-1){20}}
\put(150,20){\line(0,-1){10}}
\put(150,10){\line(1,-1){10}}

\put(143,-5){$-\frac12-\i$}
\put(170,-5){$-\i$}
\put(126,6){$-1-\frac12\i$}
\put(140,20){$-1$}
\put(167,42){$\i$}
\put(181,42){$\frac12+\i$}
\put(191,31){$1+\frac12\i$}
\put(192,18){$1$}

\end{picture}
\caption{$T$ (left), $\B$ (center)  and $U$ (right)}
\end{figure}

\begin{exa}{\rm
Let $T=T(-\i,1,1+\i,\i,-1)$. Then, its asymmetrization is $$\B=\B(1,2+2\i,\i,-1+\i)$$
and its symmetrization is 
$$U=T\left(-\frac12-\i,-\i,1,1+\frac12\i,\frac12+\i,\i,-1,-1-\frac12\i\right).$$
See Figure 3. 
}\end{exa}

Let $\beta$ be a pre-edge with $\arg(\beta)=\theta$. 
Define a random variable $\XX_\beta$ as 
$$
\XX_\beta(\omega)=\frac{|\beta|}{2}|\sin(\theta-\omega)|, 
$$
and for a pre-edge bundle $\B=\B(\beta_1,\cdots,\beta_m)$, let
$\XX_\B=\sum_{j=1}^m\XX_{\beta_j}$. 

\begin{thm}
\label{Real}
If $T$ is a realization of $\B$, then we have $X_T=\XX_\B$. 
\end{thm}
\prf
Since $T=T(\alpha_1,\cdots,\alpha_n)$ is a convex polygon, we have 
$X_T(\omega)=X_T(\omega+\pi)$ for any $\omega\in\Omega$. Hence, \\
\begin{align*}
&X_T(\omega)=\frac{X_T(\omega)+X_T(\omega+\pi)}{2}
=\sum_{j=1}^n\frac{X_{\alpha_j}(\omega)+X_{\alpha_j}(\omega+\pi)}{2}\\
&=\sum_{j=1}^n\frac{|\alpha_j|}{2}(\sin(\theta_j-\omega)_++\sin(\theta_j-\omega-\pi)_+)\\
&=\sum_{j=1}^n\frac{|\alpha_j|}{2}(\sin(\theta_j-\omega)_++(-\sin(\theta_j-\omega))_+)\\
&=\sum_{j=1}^n\frac{|\alpha_j|}{2}|\sin(\theta_j-\omega)|
=\sum_{j=1}^n\XX_{\iota(\alpha_j)}=\sum_{j=1}^m\XX_{\beta_j}=\XX_\B
\end{align*}
\prfend

\section{Minimum of deviation rate among the pre-edge bundles of fixed size}

A regular $n$-gon is denoted by $T_n~(n=2,3,\cdots)$. 
A pre-edge bundle $R_m=\B(\beta_1,\cdots,\beta_m)$ is said to be {\it regular} if 
$$|\beta_1|=\cdots=|\beta_m|~\mbox{and}~\arg(\beta_{j+1})-\arg(\beta_j)\equiv\pi/m
~({\rm mod}~\pi)~(j=1,\cdots,m).$$ 

\begin{thm}
\label{Reg}
$(i)$ It holds that 
$$
\delta(\XX_{R_m})=\nu_m~~(\mbox{see}~(3))
$$
for $m=1,2,\cdots$, and hence, 
$\delta(\XX_{R_1})<\delta(\XX_{R_2})<\delta(\XX_{R_3})<\cdots$.\\
$(ii)$~If $m\ge 2$ is even, then let $n(<m)$ be the odd number such that either 
$n=\frac{m}{2}$ or $n=\frac{m}{2}+1$. 
If $n=\frac{m}{2}$, then $\delta(X_{T_m})=\delta(X_{T_{n}})=\delta(\XX_{R_n})$ holds, 
and if $n=\frac{m}{2}+1$, then 
$\delta(X_{T_m})=\delta(\XX_{R_{m/2}})>\delta(\XX_{R_n})=\delta(X_{T_n})$ holds.  
\end{thm}
\prf\\
$(ii)$ holds since both of  $T_n$ and $T_{2n}$ are realizations of $R_n$ if 
$n=\frac{m}{2}$ is odd, and $\delta(X_{T_n})=\delta(X_{T_{2n}})=\delta(X_{T_m})
=\delta(\XX_{R_n})$ by Theorem \ref{Real}. If $n=\frac{m}{2}+1$, then by the monotonicity 
in $(i)$, $\delta(X_{T_n})=\delta(\XX_{R_n})<\delta(\XX_{R_{m/2}})$.

By Theorem \ref{Real}, to prove $(i)$, it is sufficient to prove that 
$$
\delta(X_{T_{2m}})=\nu_m~~(m=1,2,\cdots)
$$
for $T_{2m}=T(e^{\i\pi/2m},e^{\i3\pi/2m},\cdots,e^{\i\pi(1+2(2m-1))/2m})$. 
Using the diagonal representation (4), we have
\begin{align*}
&\E(X_{T_{2m}})
=\frac{2m}{2\pi}\int_{-\pi/2}^{-\pi/2+\pi/m}2\sin(\frac{\pi}{2m}-\omega)d\omega\\
&=\frac{2m}{\pi}\int_{\pi/(2m)}^{-\pi/(2m)}\sin(\frac{\pi}{2}-\omega)d\omega
=\frac{2m}{\pi}\int_{-\pi/(2m)}^{\pi/(2m)}\cos\omega d\omega
=\frac{4m}{\pi}\sin\frac{\pi}{2m}
\end{align*}
\begin{align*}
&\E(X^2_{T_{2m}})=\frac{2m}{2\pi}\int_{-\pi/2}^{-\pi/2+\pi/m}4\sin^2(\frac{\pi}{2m}-\omega)
d\omega\\
&=\frac{4m}{\pi}\int_{\pi/(2m)}^{-\pi/(2m)}\sin^2(\frac{\pi}{2}-\omega)d\omega
=\frac{4m}{\pi}\int_{-\pi/(2m)}^{\pi/(2m)}\cos^2\omega d\omega\\
&=\frac{2m}{\pi}\sin\frac{\pi}{m}+2
\end{align*}
Hence, 
\begin{align*}
&\delta(X_{T_{2m}})=\sqrt{\frac{\E(X^2_{T_{2m}})}{\E(X_{T_{2m}})^2}-1}
=\sqrt{\frac{\pi\sin\frac{\pi}{m}}{8m\sin^2\frac{\pi}{2m}}
+\frac{\pi^2}{8m^2\sin^2\frac{\pi}{2m}}-1}\\
&=\sqrt{\frac{\pi}{4m\tan(\frac{\pi}{2m})}+\frac{\pi^2}{8m^2\sin^2(\frac{\pi}{2m})}-1}.
\end{align*}

Let $x=\frac{\pi}{2m}$ and $I$ be the term inside the root in the above formula. 
Then, we have 
$$
I=\frac{x}{2\tan x}+\frac{x^2}{2\sin^2 x}-1.
$$
We show that $I$ is an increasing function of $x\in(0,\frac{\pi}{2}]$. 
We have 
\begin{align*}
&I'(x)=\frac{\cos x\sin x-x}{2\sin^2 x}+\frac{x\sin x-x^2\cos x}{\sin^3 x}\\
&=\frac{\cos x\sin^2 x+x\sin x-2x^2\cos x}{2\sin^3 x}
\end{align*}
Since 
$$
\cos x\le 1-\frac12 x^2+\frac{1}{24}x^4~~\mbox{and}~~\sin x\ge x-\frac16 x^3,
$$
\begin{align*}
&(2\sin^3x)I'(x)=\cos x\sin^2 x+x\sin x-2x^2\cos x\\
&\ge x(x-\frac16 x^3)-(2x^2-\sin^2 x)(1-\frac12 x^2+\frac{1}{24}x^4)\\
&\ge x(x-\frac16 x^3)-2x^2(1-\frac12 x^2+\frac{1}{24}x^4)
+(x-\frac16 x^3)^2(1-\frac12 x^2+\frac{1}{24}x^4)\\
&=\frac{11}{72}x^6-\frac{1}{36}x^8+\frac{1}{864}x^{10}=\frac{1}{72}x^6(12x^2-2x+11)
\end{align*}
which is positive on $x\in(0,\frac{\pi}{2}]$. 
Thus, $I(x)$ is strictly increasing in $x$, and hence, $\delta(X_{T_m})$ is strictly 
decreasing in $m=1,2,\cdots$. 

By a numerical calculation, we have $\delta(\XX_{R_m})$ as follows.
$$
\begin{array}{ccccccccc}
m&|&1&2&3&4&5&6&7\\
\hline
\delta&|&0.48342&0.09772&0.04196&0.02333&0.01485&0.01028&0.00754
\end{array}
$$
\prfend

\begin{thm}
\label{Main0}
The deviation rate $\delta(\XX_\B)$ attains the minimum among $\B\in\Xi_m$ 
if and only if $\B=R_m$. 
\end{thm}

\prf
We use the induction on $m$. If $m=1$, the statement is clear since 
$R_1$ is essentially the only element in $\Xi_1$. 

Let $m\ge2$ and assume that the statement holds for 
$\Xi_j~(j=1,\cdots,m-1)$. The boundary of the closure of $\Xi_m$ consists of 
$\cup_{j=1}^{m-1}\Xi_j$, and the $\delta$-values there are  
larger than $\delta(\XX_{R_m})$ since the assumption of the induction and 
Theorem 4. Hence, there is $\B_0\in\Xi_m$ attaining the minimum of 
$\delta(\XX_\B)$ in $\Xi_m$. We prove that $\B_0=R_m$. 

For this purpose, we take the 
symmetrization $T_0\in\Theta_{2m}$ of $\B_0$. Then, $\delta(X_{T_0})$ 
is minimum among $\delta(X_T)$ for symmetric $T\in\Theta_{2m}$. This is equivalent 
to say that $\kappa(X_{T_0})$ is minimum among $\kappa(X_T)$ for symmetric 
$T\in\Theta_{2m}$, where $\kappa(X_T)=\frac{\E(X_T^2)}{\E(X_T)^2}$.  
We'll conclude from this that $T_0$ is the regular $2m$-gon. 

Let $T_0=T(\alpha_1,\cdots,\alpha_{2m})$. 
Consider the diagonal representation (4) of $X_{T_0}$. 
Then, for any $j=1,\cdots,2m$, we have 
$$X_{T_0}(\omega)=2|\alpha_j|\sin(\arg(\alpha_j)-\omega)
~\mbox{if}~\omega\in\Omega_j,$$
where 
$$\Omega_j=\{\omega\in\Omega;
~\arg(\alpha_{j-1}-\alpha_j)<\omega\le\arg(\alpha_j-\alpha_{j+1})\}.$$
 
For a fixed $j=1,\cdots,m$ and a real number $\lambda$ near 0, let 
$$T^\lambda_0=T(\alpha_1,\cdots,(1+\lambda)\alpha_j,\alpha_{j+1},
\cdots,(1+\lambda)\alpha_{j+m},\cdots,\alpha_{2m}).$$
By the minimality, we must have 
$\frac{d}{d\lambda}\kappa(X_{T_0^\lambda})|_{\lambda=0}=0$. 

Let
$$A=\E(X_{T_0}),~B=\E(X^2_{T_0}),~a_j=\E(X_{T_0}1_{\Omega_j}),~
b_j=\E(X^2_{T_0}1_{\Omega_j}).
$$
Then, we have 
$$
\frac{\E({X_{T_0^\lambda}}^2)}{\E(X_{T_0^\lambda})^2}
=\frac{B-2b_j+2(1+\lambda)^2b_j}{(A-2a_j+2(1+\lambda)a_j)^2}+o(\lambda).
$$
Therefore, we have
$$
0=\frac{d}{d\lambda}\kappa(X_{T_0^\lambda})|_{\lambda=0}
=\frac{4b_jA-4a_jB}{A^3},
$$
and hence, 
$$\frac{b_j}{a_j}=\frac{B}{A}~~\mbox{for any}~~j=1,\cdots,2m.$$ 

Denoting 
$$
u_j=\arg(\alpha_j)-\arg(\alpha_{j-1}-\alpha_j)-\frac{\pi}{2},~
v_j=\arg(\alpha_j)-\arg(\alpha_j-\alpha_{j+1})-\frac{\pi}{2},
$$
it holds that 
\begin{align*}
&a_j=\E(X_{T_0}1_{\Omega_j})
=\int_{\arg(\alpha_{j-1}-\alpha_j)}^{\arg(\alpha_j-\alpha_{j+1})}
2|\alpha_j|\sin(\arg(\alpha_j)-\omega)\frac{d\omega}{2\pi}\\
&=\frac{|\alpha_j|}{\pi}\int_{u_j+\pi/2}^{v_j+\pi/2}\sin\omega~(-1)d\omega
=\frac{|\alpha_j|}{\pi}(\cos(v_j+\pi/2)-\cos(u_j+\pi/2))\\
&=\frac{|\alpha_j|}{\pi}(\sin u_j-\sin v_j)
\end{align*}
\begin{align*}
&b_j=\E(X_{T_0}^21_{\Omega_j})
=\int_{\arg(\alpha_{j-1}-\alpha_j)}^{\arg(\alpha_j-\alpha_{j+1})}
2|\alpha_j|^2\sin^2(\arg(\alpha_j)-\omega)\frac{d\omega}{2\pi}\\
&=\frac{|\alpha_j|^2}{\pi}\int_{u_j+\pi/2}^{v_j+\pi/2}
\sin^2(\omega)~(-1)d\omega
=\frac{|\alpha_j|^2}{\pi}\left(\frac{u_j-v_j}{2}+\frac{\sin(2v_j+\pi)-\sin(2u_j+\pi)}{4}\right)
\\&=\frac{|\alpha_j|^2}{\pi}\left(\frac{u_j-v_j}{2}+\frac{\sin 2u_j-\sin 2v_j}{4}\right). 
\end{align*}

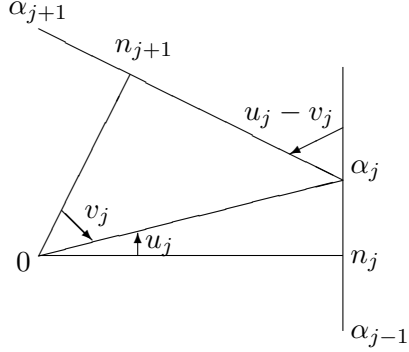
\begin{figure}[h]
\setlength{\unitlength}{1mm}
\begin{picture}(80,40)(-45,10)
\put(0,50){\line(2,-1){40}}
\put(40,45){\line(0,-1){35}}

\put(0,20){\line(1,2){12}}
\put(0,20){\line(1,0){40}}
\put(0,20){\line(4,1){40}}

\put(-3,18){0}
\put(41,31){$\alpha_j$}
\put(41,9){$\alpha_{j-1}$}
\put(41,19){$n_j$}
\put(-4,52){$\alpha_{j+1}$}
\put(10,47){$n_{j+1}$}

\put(3,26){\vector(1,-1){4}}
\put(6,25){$v_j$}
\put(13,20){\vector(0,1){3}}
\put(14,21){$u_j$}
\put(40,37){\vector(-2,-1){7}}
\put(27,38){$u_j-v_j$}

\end{picture}
\caption{$\alpha_{j-1},\alpha_{j},\alpha_{j+1},n_j,n_{j+1},u_j,v_j$}
\end{figure}

For $j=1,\cdots,2m$, let
$$
n_j=\mbox{the perpendicular leg from 0 to the line}~\alpha_{j-1}\alpha_j
$$
Then, the above $u_j$ and $v_j$ have another representation (see Figure 6) that
\begin{align}
&u_j=\arg(\alpha_j)-\arg(n_j),~v_j=\arg(\alpha_j)-\arg(n_{j+1})\in(-\pi,\pi]\\
&u_j-v_j=\mbox{the exterior angle at}~\alpha_j>0\nonumber
\end{align} 
and we have 
\begin{equation}
\frac{B}{A}=\frac{|\alpha_j|(2u_j-2v_j+\sin2u_j-\sin2v_j)}{4(\sin u_j-\sin v_j)}
~\mbox{for}~j=1,\cdots,2m.
\end{equation} 

For a fixed $j=1,\cdots,m$ and a real number $\lambda$ near 0, let 
$$T^{\i\lambda}_0=T(\alpha_1,\cdots,(1+\i\lambda)\alpha_j,\alpha_{j+1},\cdots,
(1+\i\lambda)\alpha_{j+m},\cdots,\alpha_{2m}).$$
Then, we have 
$$
\frac{\E({X_{T_0^\lambda}}^2)}{\E(X_{T_0^\lambda})^2}
=\frac{B+2d_j^2\frac{\lambda}{2\pi}-2d_{j+1}^2\frac{\lambda}{2\pi}}
{(A+2d_j\frac{\lambda}{2\pi}-2d_{j+1}\frac{\lambda}{2\pi})^2}+o(\lambda),
$$
where 
\begin{align*}
&d_j=X_{T_0}(\arg(\alpha_{j-1}-\alpha_j))\\&~~~=2|\alpha_j|
\sin(\arg(\alpha_j)-\arg(\alpha_{j-1}-\alpha_j))=2|\alpha_j|\cos u_j\\
&d_{j+1}=X_{T_0}(\arg(\alpha_j-\alpha_{j+1}))\\&~~~~~~
=2|\alpha_j|\sin(\arg(\alpha_j)-\arg(\alpha_j-\alpha_{j+1}))=2|\alpha_j|\cos v_j.
\end{align*}
Therefore, we have
$$
0=\frac{d}{d\lambda}\kappa(X_{T_0^{\i\lambda}})|_{\lambda=0}
=\frac{\frac{1}{\pi}(d_j^2-d_{j+1}^2)A-\frac{2}{\pi}(d_j-d_{j+1})B}{A^3},
$$
and hence, either $d_j=d_{j+1}$ or $\frac{d_j+d_{j+1}}{2}=\frac{B}{A}$. 

If $\frac{d_j+d_{j+1}}{2}=\frac{B}{A}$ holds, then 
by (6), we have 
$$
\frac{|\alpha_j|(2u_j-2v_j+\sin2u_j-\sin2v_j)}{4(\sin u_j-\sin v_j)}
=|\alpha_j|(\cos u_j+\cos v_j).
$$
Hence,
$$
2(u_j-v_j)=\sin2u_j-\sin2v_j+4\sin(u_j-v_j).
$$
Since $u_j-v_j>0$ and if $u_j-v_j\le\frac{\pi}{2}$, then $\sin(u_j-v_j)\ge\frac{2}{\pi}(u_j-v_j)$, 
we have a contradiction that $\sin2u_j-\sin2v_j+4\sin(u_j-v_j)>2(u_j-v_j)$. 
Therefore, $u_j-v_j>\frac{\pi}{2}$ and  
the exterior angle of $\alpha_j$ is larger than $\pi/2$. By the symmetry, 
the exterior angle of $\alpha_k$ for $k\ne j,j+\pi$ is smaller than $\pi/2$. 
Therefore, for these $k$, $\frac{d_j+d_{j+1}}{2}=\frac{B}{A}$ is impossible, and hence, 
$d_k=d_{k+1}$ holds. Thus, $u_k=-v_k>0$. This implies $|n_k|=|n_{k+1}|$ 
for any $k\ne j,j+m$. 

Any case, we have $|n_k|=|n_{k+1}|$ except for $k=j,j+m$. 
This implies that there are 2 classes 
\begin{align*}
&|n_{j+1}|=|n_{j+2}|=\cdots=|n_{j+m}|\\
&|n_{j+m+1}|=|n_{j+m+2}|=\cdots=|n_{j+2m}|\\
&(\mbox{suffices are considered modulo}~2m).
\end{align*}
By symmetry, the values of these 2 classes coincide. 

Thus, $|n_1|=|n_2|=\cdots=|n_{2m}|$, which implies that 
$T_0$ has a inscribed circle with radius $r:=|n_1|$. Hence, 
$u_j=-v_j~(j=1,\cdots,2m)$. 
Then by (6), 
\begin{align*}
&\frac{B}{A}=\frac{|\alpha_j|(2u_j+\sin 2u_j)}{4\sin u_j}
=\frac{|\alpha_j|\cos u_j(2u_j+\sin 2u_j)}{4\sin u_j\cos u_j}\\
&=\frac{r(2u_j+\sin 2u_j)}{2\sin 2u_j}=\frac{r}{2}~\frac{2u_j}{\sin 2u_j}+\frac{r}{2}~~
(j=1,\cdots,2m).
\end{align*}
It follows from this that $\frac{2u_j}{\sin 2u_j}$ is the same for $j=1,\cdots,2m$. 
Since the correspondence $x\mapsto\frac{x}{\sin x}$ for $x>0$ is one-to-one, 
we have $u_1=\cdots=u_{2m}$. Thus, $T_0$ has the same exterior angle $2u_j$ 
at vertex $\alpha_j$ for $j=1,\cdots,2m$. Together with that $T_0$ has a inscribed circle, 
$T_0$ is a regular $2m$-gon, which completes the proof. 
\prfend 

\section{Minimum of deviation rate among $n$-gons}

\begin{lem}
\label{Cyc}
If $n\ge 2$, then the set
$$
P(n):=\left\{ (c_0,\dots,c_{n-1}) \in \{-1,1\}^n \ \left|\ 
\sum_{j=0}^{n-1} c_j\exp\left(\frac{j\pi\i}{n}\right)=0
\right.\right\}
$$
is empty if and only if $n$ is a power of $2$.
\end{lem}

\prf
If $n$ is not a power of $2$, then take an odd factor $p$ of $n$. 
Denote $j=0,1,\cdots,n-1$ as
$$j=k\frac{n}p+\ell~~(\ell=0,1,\cdots,\frac{n}{p}-1;~k=0,1,\cdots,p-1).$$
Given $c_\ell\in\{-1,1\}~(\ell=0,1,\cdots,\frac{n}{p}-1)$ arbitrary, define 
$$c_j=c_\ell(-1)^k~~(j=0,1,\cdots,n-1).$$ 
Then, we have
\begin{align*}
&\sum_{j=0}^{n-1} c_je^{\frac{j\pi\i}{n}}
=\sum_{j=0}^{n-1}c_\ell(-1)^k\exp\left(\frac{(k\frac{n}{p}+\ell)\pi\i}{n}\right)\\
&=\sum_{\ell=0}^{(n/p)-1}c_\ell\exp\left(\frac{\ell\pi\i}{n}\right)
\sum_{k=0}^{p-1}(-1)^k\exp\left(\frac{k\pi\i}{p}\right)\\
&=\sum_{\ell=0}^{(n/p)-1}c_\ell\exp\left(\frac{\ell\pi\i}{n}\right)
\left(\sum_{k=0}^{\frac{p-1}{2}}\exp\left(\frac{(2k\i)\pi}{p}\right)
-\sum_{k=1}^{\frac{p-1}{2}}\exp\left(\frac{(2k-1)\pi\i}{p}\right)\right)\\
&=\sum_{\ell=0}^{(n/p)-1}c_\ell\exp\left(\frac{2\pi\ell\i}{n}\right)
\left(\sum_{k=0}^{\frac{p-1}{2}}\exp\left(\frac{2k\pi\i}{p}\right)
+\sum_{k=1}^{\frac{p-1}{2}}\exp\left(\frac{(2k-1+p)\pi\i}{p}\right)\right)
\end{align*}
\begin{align*}
&=\sum_{\ell=0}^{(n/p)-1}c_\ell\exp\left(\frac{2\pi\ell\i}{n}\right)
\left(\sum_{k=0}^{(p-1)/2}\exp\left(\frac{2k\pi\i}{p}\right)
+\sum_{k=(p+1)/2}^{p-1}\exp\left(\frac{2k\pi\i}{p}\right)\right)\\
&=\sum_{\ell=0}^{(n/p)-1}c_\ell\exp\left(\frac{2\pi\ell\i}{n}\right)
\sum_{k=0}^{p-1}\exp\left(\frac{2k\pi\i}{p}\right)=0.
\end{align*}
Therefore $P(n)$ contains a non-empty subset
$$
Q(p):= \left\{(c_0,\dots,c_{n-1})\in \{-1,1\}^n;~c_{k(n/p)+\ell}= c_{\ell} (-1)^k\right\}
$$
with $\mathrm{Card}(Q(p))=2^{n/p}$. Thus we see that $P(n)$ is non-empty. 

Next assume that $n=2^s$ with $s\ge 1$. Since 
$x^n+1\in\Z[x]$ is the minimum polynomial of $w:=\exp\left(\i\frac{\pi}{n}\right)$, 
$$
\sum_{j=0}^{n-1} c_j\exp\left(\frac{j\pi\i}{n}\right)=c_0+c_1w+\cdots+c_{n-1}w^{n-1}
$$
cannot be 0 for any $(c_0,c_1,\cdots,c_{n-1})\in\{-1,1\}^n$. 
Hence, $P(n)=\emptyset$.
\prfend

\begin{thm}
\label{Main}
Let $n\ge 2$ be an integer which is not a power of $2$. Then we have
$$
\min_{T\in\Theta_n}\delta(X_T)=\delta(\XX_{R_n})=\nu_n~~(\mbox{see}~(3)).
$$ 
The minimum is attained if and only if the asymmetrization of $T$ is $R_n$, 
and hence, if and only if $T$ is similar to the
polygon $T(\alpha_1,\dots, \alpha_n)$ with 
$$
\{\alpha_{j+1}-\alpha_j;~j=1,\dots,n\}
=\left\{c_j\exp\left(\frac{j\pi\i}{n}\right);~(c_0,\dots,c_{n-1})\in P(n)\right\}.
$$
Here $P(n)$ is defined in Lemma \ref{Cyc}.
\end{thm}

\prf
In light of Theorem \ref{Main0}, 
$T^*\in\Theta_n$ attains the minimum of $\delta(X_T)$
among $T\in\Theta_n$ if $T^*$ is a realization of the pre-edge bundle $R_n$.
A realization of $R_n$, say $T$, may have more than $n$ edges 
but this number is reduced to $n$ if and only if there exists a one-to-one correspondence between the set of pre-edges $\{\beta_1,\cdots,\beta_n\}$ 
of $R_n$ and the set of edges of $T$, 
i.e., there exists 
$(c_0,c_1,\dots, c_{n-1})\in \{-1,1\}^n$ such that 
$\{c_0\beta_0,c_1\beta_1,\dots, c_{n-1}\beta_{n-1}\}$ is the set of edges of $T$.
Since we may assume $\beta_j=\exp\left(\frac{j\pi\i}{n}\right)$ for $i=0,\dots, n-1$,
this condition is satisfied if and only if
\begin{equation}
\label{Cond}
    \sum_{j=0}^{n-1} c_j\exp\left(\frac{j\pi\i}{n}\right)=0
\end{equation}
is solvable in $c_j\in \{-1,1\}$, i.e. $P(n)$ is non-empty.
\prfend

%





Hereafter $n>1$ is always assumed to be an integer which is not a power of $2$.
There is a natural map $\sigma$ and $\tau$ from $P(n)$ to itself defined by
$$
\sigma((c_0,c_1,\dots, c_{n-1}))=(c_1,\dots c_{n-1}, -c_0)
$$
and
$$
\tau((c_0,c_1,\dots, c_{n-1}))=(c_{n-1},c_{n-2},\dots, c_0)
$$
which corresponds to the symmetry of dihedral group $D_{2n}$:
the rotation of angle $\pi/n$ and the reflexion. 
Two elements $(c_0,\dots,c_{n-1})$ and $(c'_0,\dots, c'_{n-1})$
of $P(n)$ give congruent realizations if and only if
$$
(c_0,\dots, c_{n-1}) = \sigma^j  (c'_0,\dots, c'_{n-1})
$$
or
$$
(c_0,\dots, c_{n-1}) = \sigma^j \tau (c'_0,\dots, c'_{n-1})
$$
for some $j\in \{0,\dots, 2n-1\}$. 
Lemma \ref{Cyc} and Theorem \ref{Main} can be restated in a geometric form.

\begin{thm}
\label{MainGeom}
If $1<n\in \N$ is not a power of $2$, then the minimum polygon in $\Theta_n$ is a Reinhardt $n$-gon and vice versa.
A Reinhardt $n$-gon exists if and only if $n$ is not a power of $2$.
\end{thm}

Though it is not stated in this manner, the latter statement follows from the characterization of Reinhardt \cite{RH}, see the discussion after the proof. 
\bigskip

\prf
Assume that $1<n\in \N$ is not a power of $2$. 
By Theorem \ref{Main}, the minimum $n$-polygon $T^*$ can be considered, without 
generality, to have the asymmetrization 
$$
\B\left(c\exp\left(\frac{\pi\i}{2n}\right),c\exp\left(\frac{3\pi\i}{2n}\right),
\cdots,c\exp\left(\frac{(2n-1)\pi\i}{2n}\right)\right)
$$
for some appropriate $c>0$, so that its symmetrization is
$$
T_{2n}=T\left(e^{\frac{(-n)\pi\i}{2n}},e^{\frac{(-n+2)\pi\i}{2n}},\cdots,
e^{\frac{(3n-2)\pi\i}{2n}}\right). 
$$
Let $T^*=T(P_1,P_2,\cdots,P_n)$. 
Since $X_{T^*}=X_{T_{2n}}$ and by (4), $X_{T^*}(\omega)$ is a periodic function of period $\pi/n$ such that if $\omega\equiv\eta~(\mbox{mod}~\pi/n)$ with 
$\eta\in(-\frac{\pi}{2n},\frac{\pi}{2n}]$,  
$$
X_{T^*}(\omega)=2\cos \eta.
$$
Therefore, $X_{T^*}(\omega)$ for $\omega\in(0,2\pi]$ repeats its maximum 2 and 
its minimum $2\cos\left(\frac{\pi}{2n}\right)$, $2n$ times. 

Since $T^*$ has the asymmetrization of the same size $n$, there are no parallel edges 
in $T^*$. Hence, the endpoints of the minimum shadow of $T^*$, 
say at $\omega=\frac{j\pi}{n}+\frac{\pi}{2n}$ come from a vertex, say $P_k$, 
and an edge, say $P_lP_{l+1}$. Those of the 2 neighboring maximum shadows come 
from the vertices $P_k,P_l$ and from the vertices $P_k,P_{l+1}$, respectively. 
Hence, $\overline{P_kP_l}=\overline{P_kP_{l+1}}$ 
holds. Replace the edge $P_lP_{l+1}$ by the circular arc centered at $P_k$ 
having the endpoints at $P_l$ and $P_{l+1}$. Repeating this replacement 
for $j=1,2,\cdots,n$, we get a convex body of constant width. 
It is easy to see that this convex body is a Reuleaux $p$-body for 
some $p$. Hence, $T^*$ is a Reinhardt $n$-gon. 
 
We prove the converse. 
Let $T$ be a Reinhardt $n$-gon coming from a Reuleaux 
$p$-body of width $r$. 
By (2), we can rearrange the circular arc of $S$ by parallel translations  
into a circle $\S$ of radius $r$ so that the circular arcs are essentially 
disjoint and cover just a half part of $\S$, and by the rotation of 
angle $\pi$, they moved to the other half part of $\S$. The edges of $T$ 
correspond to the chord of the circular arcs. By the above property, 
it is easy to see that the asymmetrization of $T$ is a regular 
pre-edge bundle of size $n$. Hence, $\delta(X_T)=\nu_n$ and $T$ is 
the minimum polygon. 

The ``if'' part of the last statement follows from the first part. 
To prove the ``only if'' part, 
suppose that a Reinhardt $n$-polygon $T$ exists for $n=2^k$. 
Then by the above argument, $T$ has the asymmetrization $R_n$ which 
has a realization $T$ of the same size. This contradicts Lemma 2. 
\prfend

Let us describe the correspondence between $(c_0,\dots, c_{n-1})\in P(n)$
and cyclic integer vector. 
Choose $i\in \{0,1,\dots, 2n-1\}$ that
$\sigma^i(c_0,\dots, c_{n-1})=(d_0,\dots, d_{n-1})$ with $d_0=d_{n-1}$.
Count number of runs of 
$1$ and $-1$ in $(d_0,\dots,d_{n-1})$, i.e., 
we write $(d_0,\dots,d_{n-1})$ like $1^{n_1}(-1)^{n_2}\dots$ or $(-1)^{n_1}1^{n_2}\dots$. Then $(n_1,n_2,\dots, n_{p})$ is the desired cyclic vector.
For the converse direction, we choose either
$1^{n_1}(-1)^{n_{2}}\dots$ or $(-1)^{n_1}1^{n_{2}}\dots$.

Reinhardt \cite{RH} gave an alternative
characterization of the cyclic vector $(n_1,n_2,\dots, n_{p})$ with $n=\sum_{i=1}^{p} n_i$: 
it is a cyclic vector if and only if
$p$ is odd and the polynomial
$$
    1-z^{n_1}+z^{n_1+n_2}-\dots + z^{n_1+n_2+\dots +n_{p-1}}
$$
is divisible by $\Phi_{2n}(z)$, the $2n$-th cyclotomic polynomial. For completeness, we show that this characterization is equivalent to ours.
As above, we assume that $d_0=d_{n-1}$. 
From (\ref{Cond}) we have
$$
 \sum_{j=0}^{n-1} d_j z^j=0\quad \text{with } z=\exp\left(\frac{\pi\i}{n}\right).
$$
We consider that $z$ is a variable and 
multiply $z-1$. Then we see
\begin{align}
&(z-1)\sum_{j=0}^{n-1} d_j z^j\notag\\
&=1-2z^{n_1}+2z^{n_1+n_2}-\dots + 2z^{n_1+\dots+n_{p-1}}-z^{n} \label{RHeq}\\
&\equiv
2-2z^{n_1}+2z^{n_1+n_2}-\dots + 2z^{n_1+\dots+n_{p-1}}
\pmod{\Phi_{2n}(z)}\notag.
\end{align}
Dividing by $2$, we see the condition of Reinhardt. To get the converse, we just go backwards. Note that the polynomial (\ref{RHeq})
is divisible by $z-1$ as $p$ is odd. By Lemma \ref{Cyc}, the second statement of Theorem \ref{MainGeom} is derived from the Reinhardt criterion.
The subset $Q(p)$ in the proof of Lemma \ref{Cyc}
corresponds to $p$-fold rotational symmetry. We can find 
Reinhardt polygons without 
any symmetry \cite{HM13,HM19}.

\section{Truncation of the regular triangle}

Let $X,Y$ be general $\R$-valued, square integrable random variables on the 
probability space $\Omega$. Assume further that $X\ge0$ everywhere and 
$\E(X)>0$. Recall that 
$$\kappa(X)=\frac{\E(X^2)}{\E(X)^2}=\delta(X)^2+1.$$

It holds for any $t\in\R$ with sufficiently small modulus that 
\begin{align*}
&\kappa(X+tY)=\frac{\E((X+tY)^2)}{\E(X+tY)^2}
=\frac{\E(X^2)}{\E(X)^2}\frac{1+2t\frac{\E(XY)}{\E(X^2)}+t^2\frac{\E(Y^2)}{\E(X^2)}}
{1+2t\frac{\E(Y)}{\E(X)}+t^2\frac{\E(Y)^2}{\E(X)^2}}\\
&=\frac{\E(X^2)}{\E(X)^2}\left(1+2t\left(\frac{\E(XY)}{\E(X^2)}
-\frac{\E(Y)}{\E(X)}\right)+O(t^2)\right).
\end{align*}
If $\frac{\E(XY)}{\E(X^2)}-\frac{\E(Y)}{\E(X)}=0$, then we have
$$
\kappa(X+tY)=\frac{\E(X^2)}{\E(X)^2}\left(1+t^2\left(\frac{\E(Y^2)}{\E(X^2)}
-\frac{\E(Y)^2}{\E(X)^2}\right)+O(t^3)\right). 
$$
Hence, the following Lemma holds.
\begin{lem}
\label{L2}
$(i)$ $\frac{d}{dt}\delta(X+tY)|_{t=0}>,=,<0$ if and only if 
$\frac{\E(XY)}{\E(X^2)}>,=,<\frac{\E(Y)}{\E(X)}$, respectively. \\
$(ii)$ Assume that ``$=$'' holds in $(i)$. Then, there exists $\varepsilon>0$ 
such that 
$\delta(X+tY)>,<\delta(X)$ for any $t\in(-\varepsilon,\varepsilon)\setminus\{0\}$ if 
$\frac{\E(Y^2)}{\E(X^2)}>,<\frac{\E(Y)^2}{\E(X)^2}$, respectively.
\end{lem}
\prf
$(i)$ follows since
$$
\textstyle\frac{d}{dt}\delta(X+tY)=(1/2)(\kappa(X+tY)-1)^{-1/2}\frac{d}{dt}\kappa(X+tY).
$$
$(ii)$ follows since $\delta(X+tY)$ is a monotone increasing function of $\kappa(X+tY)$.
\prfend

\begin{lem}
\label{L3}
Let $T$ and $S$ be triangles in $\C$. Then, $X_T=X_S$~(a.s.)  
holds if and only if there exists $z\in\C$ such that either $S=T+z$ or 
$S=-T+z$. 
\end{lem}
\prf
$X_{-T}=X_T$ holds since they have the same asymmetrization. Hence, the ``if'' part holds. 

Let us prove the ``only if'' part. Let $T=T(\alpha,\beta,\gamma)$. 
Consider $X_T(\omega)$ as a function of $\omega\in\R/\Z$. 
Then, it is locally minimal if and only if $\omega$ is equal to either of 
$$\pm(\arg(\beta)-\arg(\alpha)),~\pm(\arg(\gamma)-\arg(\beta)),
~\pm(\arg(\alpha)-\arg(\gamma))$$
modulo $2\pi$. If $X_T=X_S$~(a.s.), then they should have the same set of $\omega$ 
as this. Also, at any of these $\omega$, they should have the same height. This implies 
that either $S=T+z$ or $S=-T+z$ for some $z\in\C$. 
\prfend

\begin{lem}
\label{L4}
Let $T$ and $S$ be triangles in $\C$. Then, 
$$\E(X_TX_S)\le\E(X_T^2)^{1/2}\E(X_S^2)^{1/2}.$$ 
The equality holds if and only if there exist $\lambda\ge0$ and $z\in\C$ such that 
either $S=\lambda T+z$ or $S=-\lambda T+z$. 
\end{lem}
\prf
This is the Cauchy-Schwartz inequality for
the inner product $\langle X,Y \rangle=E(XY)$.
The equality holds if and only if there exists 
$\lambda>0$ such that $X_S=\lambda X_T=X_{\lambda T}$~(a.s.).  
Hence, by Lemma \ref{L3}, if and only if $S=\lambda T+z$ or $S=-\lambda T+z$ 
for some $z\in\C$.
\prfend

\begin{lem}
\label{L5}
Let $T=T(\alpha,\beta,\gamma)$ be a regular triangle. 
Let $\sigma\in\C\setminus\{0\}$. Then, 
$$
\frac{\E(X_TX_\sigma)}{\E(X_T^2)}\le\frac{\E(X_\sigma)}{\E(X_T)}.
$$ 
The equality holds if and only if $\vec{0\sigma}$ is parallel to
one of the edges of $T$, that is, $\arg(\sigma)$ is equal to one of  
$$\pm(\arg(\beta)-\arg(\alpha)),~\pm(\arg(\gamma)-\arg(\beta)),
~\pm(\arg(\alpha)-\arg(\gamma))$$
modulo $2\pi$. 
\end{lem}
\prf
Without loss of generality, we assume that the length of the edges of $T$ is 1. 
It holds that $\E(X_{\mu z_1}X_{\mu z_2})=\E(X_{z_1}X_{z_2})$ for any 
$\mu, z_1,z_2\in\C$ with $|\mu|=1$. 
Moreover, since $e^{2\pi \i/3}T=T$, we have 
$$
\E(X_TX_\sigma)=\E(X_TX_{e^{2\pi \i/3}\sigma})=\E(X_TX_{e^{4\pi \i/3}\sigma}).
$$
Hence, $\E(X_TX_\sigma)=(1/3)\E(X_TX_S)$ 
with the regular triangle $S=T(0,\sigma,e^{\pi \i/3}\sigma)$.
Therefore by Lemma \ref{L4}, 
\begin{align*}
&\E(X_TX_\sigma)=(1/3)\E(X_TX_S)\le(1/3)\E(X_T^2)^{1/2}\E(X_S^2)^{1/2}\\
&=(1/3)\E(X_T^2)^{1/2}\E((|\sigma|X_T)^2)^{1/2}=(|\sigma|/3)\E(X_T^2).
\end{align*}
The equality holds if and only if $\vec{0\sigma}$ is parallel to
one of the edges of $T$. Thus, we have 
$$\frac{\E(X_TX_\sigma)}{\E(X_T^2)}\le\frac{|\sigma|}{3}=\frac{\E(X_\sigma)}{\E(X_T)}$$ 
with the equality if and only if $\vec{0\sigma}$ is parallel to one of the edges of $T$. 
\prfend

\begin{thm}
\label{Trun}
A parallel truncation of the regular triangle increases the deviation rate, while 
a non-parallel truncation decreases it. 
\end{thm}
\noindent\underline{Parallel truncation:}

Let $T=T(\alpha,\beta,\gamma)$ be a regular triangle of the edge length 1. 
We also assume that it is of the counter clockwise order. 
Let $t>0$ be sufficiently small. Let 
$$\beta_t=(1-t)\alpha+t\beta,~\gamma_t=(1-t)\alpha+t\gamma.$$
Let
$$
V_t=T(\beta_t,\beta,\gamma,\gamma_t)
$$ 
be a parallel truncation of $T$ at $\alpha$. 
Then, we have 
\begin{align*}
&X_{V_t}=X_{\beta-\beta_t}+X_{\gamma-\beta}+X_{\gamma_t-\gamma}
+X_{\beta_t-\gamma_t}\\
&=X_T-tX_{\beta-\alpha}-tX_{\alpha-\gamma}+tX_{\beta-\gamma}\\
&=X_T-tX_Y
\end{align*}
with $Y=X_{\c}+X_{\b}-X_{-\a}$, 
where we denote $\a=\gamma-\beta,~\b=\alpha-\gamma,~\c=\beta-\alpha$. 
Since $\E(X_{\c})=\E(X_{\b})=\E(X_{-\a})=1/\pi$, we have by Lemma \ref{L5}, 
$$
\frac{\E(X_TX_{\c})}{\E(X_T^2)}=
\frac{\E(X_TX_{\b})}{\E(X_T^2)}=
\frac{\E(X_TX_{-\a})}{\E(X_T^2)}=
\frac{1/\pi}{\E(X_T)},
$$
and hence, 
$$
\frac{\E(X_TY)}{\E(X_T^2)}=\frac{\E(X_TX_{\c})+\E(X_TX_{\b})-\E(X_TX_{-\a})}
{\E(X_T^2)+\E(X_T^2)-\E(X_T^2)}=\frac{1/\pi}{\E(X_T)}=\frac{\E(Y)}{\E(X_T)}.
$$

Now, we prove that 
$$
\frac{\E(Y^2)}{\E(X_T^2)}>\frac{\E(Y)^2}{\E(X_T)^2}
$$
so that $\delta(V_t)>\delta(T)$ for sufficiently small $|t|\ne0$ by Lemma \ref{L2}. 
We have
\begin{align*}
&\E(Y)=\E(X_{\b}+X_{\c}-X_{-\a})=\frac{1}{\pi}\\&\E(X_T)=\frac{3}{\pi}\\
&\E(Y^2)=\E((X_{\b}+X_{\c}-X_{-\a})^2)\\
&=\frac34+(2/4\pi)(V(2\pi/3)-2V(\pi/3))
=\frac13-\frac{\sqrt{3}}{4\pi}\\
&\E(X_T^2)=\frac34+(6/4\pi)V(2\pi/3)=\frac12+\frac{3\sqrt{3}}{4\pi}
\end{align*}
$$
\frac{\E(Y^2)}{\E(X_T^2)}=\frac{\frac13-\frac{\sqrt{3}}{4\pi}}{\frac12+\frac{3\sqrt{3}}{4\pi}}=0.214\cdots>\frac19=\frac{\E(Y)^2}{\E(X_T)^2},
$$
and complete the proof that $\delta(V_t)>\delta(T)$ for sufficiently small $|t|\ne0$ 
in the case of parallel truncation.\\

\noindent\underline{Non-parallel truncation:}

Let $T=T(\alpha,\beta,\gamma)$ be a regular triangle of the edge length 1. 
We also assume that it is of the counter clockwise order. 
Let $0<\lambda\ne1$ and $t>0$ be sufficiently small. Let 
$$\beta_t=(1-t)\alpha+t\beta,~\gamma_t=(1-\lambda t)\alpha+\lambda t\gamma.$$
Let
$$
V_t=T(\beta_t,\beta,\gamma,\gamma_t)
$$ 
be a non-parallel truncation of $T$ at $\alpha$. 
Then, we have 
\begin{align*}
&X_{V_t}=X_{\beta-\beta_t}+X_{\gamma-\beta}+X_{\gamma_t-\gamma}
+X_{\beta_t-\gamma_t}\\
&=X_T-tX_{\c}-tX_{\lambda\b}+tX_{\e}\\
&=X_T-tY,
\end{align*}
where we denote $\b=\alpha-\gamma,~\c=\beta-\alpha,~\e=\c+\lambda\b$ 
and $Y=X_{\c}+X_{\lambda\b}-X_{\e}$. 

By Lemma \ref{L5}, we have
$$
\frac{\E(X_TX_{\c})}{\E(X_T^2)}=\frac{\E(X_TX_{\b})}{\E(X_T^2)}
=\frac{1/\pi}{\E(X_T)}~\mbox{and}~
\frac{\E(X_TX_{\e})}{\E(X_T^2)}<\frac{|\e|/\pi}{\E(X_T)}.
$$
Hence,
$$
\frac{\E(X_TX_{\c})+\lambda\E(X_TX_{\b})-\E(X_TX_{\e}))}
{\E(X_T^2)+\lambda\E(X_T^2)-|\e|\E(X_T^2)}>\frac{1/\pi}{\E(X_T)},
$$
and we have
$$
\frac{\E(X_TY)}{\E(X_T^2)}>\frac{(1+\lambda-|\e|)/\pi}{\E(X_T)}.
$$
Thus, by Lemma \ref{L2}, 
$\frac{d}{dt}\delta(X_T+tY)|_{t=0}>0$, which implies that 
$$\delta(X_{V_t})=\delta(X_T-tY)<\delta(X_T)$$ if $t>0$ is small, 
which completes the proof. 

\section{Remaining problems}
When $n$ is a power of $2$, the method in this paper 
does not apply because by Lemma \ref{Cyc}, 
there is no $n$-gon realization of the regular pre-edge bundle $R_n$. 
The case $n=4$ may be of special interest.
\begin{figure}[h]
\begin{center}
\includegraphics[12cm,4cm]{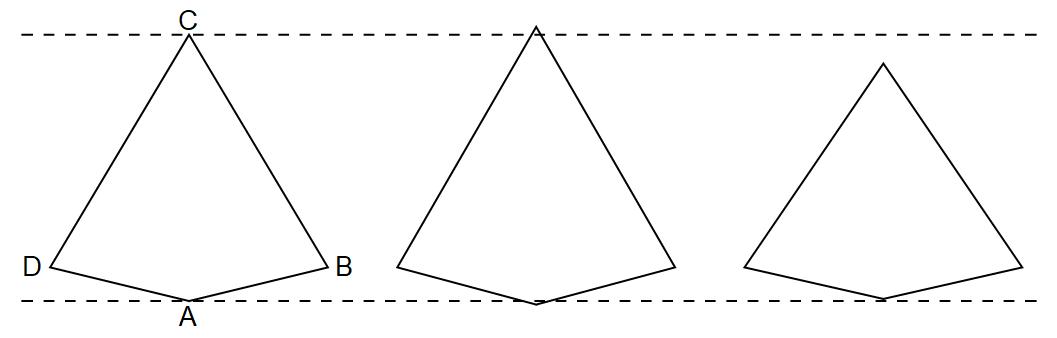}
\vspace{-1em}
\end{center}
\caption{Left: $\delta$-minimum,\quad Center: Problem 1 \& 2,\quad Right: Problem 3\label{Tetra}}
\end{figure}
By numerical calculation, a possible
minimum of $\delta$ in $\Theta_4$ is attained by the 
kite-shape with vertices 
$$(A, B, C, D)\approx ((0,-0.24213332485),(1,0), (0,1.67502597318),
(-1,0))$$ 
having the deviation rate $0.035306425305$.
This is not a solution to any of three optimization problems for $n=4$
in the introduction. 
See Figure \ref{Tetra} for a comparison of solutions having the same horizontal width.
Indeed from $AB<AC<BC<BD$, we see that this shape is not optimal for Problems 1 and 2. The optimizers for these problems are the same, which is 
based on the regular triangle with an additional vertex similar to the 
construction of 
Reinhardt polygon, see \cite{Mo, BF}.  The kite $ABCD$ is not optimal for Problem 3 either. 
Indeed \cite{AHM}
showed that the maximum width of quadrangles with the unit perimeter is
close to $\frac{\sqrt{-9 + 6\sqrt{3}}}4\approx 0.295$, possibly attained by
$$
\left(\left(0, -\frac {\sqrt{-3 + 2 \sqrt{3}}}3\right), (1, 0), \left(0, 
\sqrt{1 + \frac 2{\sqrt{3}}}\right), (-1, 0)\right),
$$
while the kite shape gives the value 0.288. For Problem 1,2 and 3, the solutions have algebraic expressions. We do not know if our kite has such an algebraic expression.

In the next paper, we will discuss the minimality of $\delta$-values and
the minimal shapes under all infinitesimal deformations
including the cases when $n$ is a power of $2$.

\end{document}